\documentclass[10pt]{amsart}

\usepackage{amsmath,amsthm,amsfonts,amscd,amssymb}
\usepackage{hyperref,wasysym}
\usepackage{verbatim}

\newtheorem{thm}{Theorem}[section]

\newtheorem{lem}[thm]{Lemma}
\newtheorem{prop}[thm]{Proposition}

\theoremstyle{definition}

\theoremstyle{remark}
\newtheorem{rem}{Remark}[section]

\begin{document}

\title[Normal bases of small height]{Normal bases of small height in Galois number fields}

\author{Lenny Fukshansky}
\author{Sehun Jeong}

\address{Department of Mathematics, 850 Columbia Avenue, Claremont McKenna College, Claremont, CA 91711}
\email{lenny@cmc.edu}
\address{Department of Mathematics, 850 Columbia Avenue, Claremont McKenna College, Claremont, CA 91711}
\email{Sehun.Jeong@ClaremontMcKenna.edu}


\begin{abstract} 
Let $K$ be a number field of degree $d$ so that $K/\mathbb Q$ is a Galois extension. The {\it normal basis theorem} states that $K$ has a $\mathbb Q$-basis consisting of algebraic conjugates, in fact $K$ contains infinitely many such bases. We prove an effective version of this theorem, obtaining a normal basis for $K/\mathbb Q$ of bounded Weil height with an explicit bound in terms of the degree and discriminant of $K$. In the case when $d$ is prime, we obtain a particularly good bound using a different method.
\end{abstract}

\maketitle

\def\A{{\mathcal A}}
\def\B{{\mathcal B}}
\def\C{{\mathcal C}}
\def\D{{\mathcal D}}
\def\F{{\mathcal F}}
\def\x{{\mathcal H}}
\def\I{{\mathcal I}}
\def\J{{\mathcal J}}
\def\K{{\mathcal K}}
\def\L{{\mathcal L}}
\def\M{{\mathcal M}}
\def\N{{\mathcal N}}
\def\O{{\mathcal O}}
\def\R{{\mathcal R}}
\def\s{{\mathcal S}}
\def\V{{\mathcal V}}
\def\W{{\mathcal W}}
\def\X{{\mathcal X}}
\def\Y{{\mathcal Y}}
\def\H{{\mathcal H}}
\def\Z{{\mathcal Z}}
\def\OO{{\mathcal O}}
\def\BB{{\mathbb B}}
\def\cee{{\mathbb C}}
\def\EE{{\mathbb E}}
\def\Nn{{\mathbb N}}
\def\pee{{\mathbb P}}
\def\que{{\mathbb Q}}
\def\real{{\mathbb R}}
\def\zed{{\mathbb Z}}
\def\hyp{{\mathbb H}}
\def\aa{{\mathfrak a}}
\def\HH{{\mathfrak H}}
\def\qbar{{\overline{\mathbb Q}}}
\def\eps{{\varepsilon}}
\def\ahat{{\hat \alpha}}
\def\bhat{{\hat \beta}}
\def\gt{{\tilde \gamma}}
\def\h{{\tfrac12}}
\def\be{{\boldsymbol e}}
\def\bei{{\boldsymbol e_i}}
\def\bff{{\boldsymbol f}}
\def\ba{{\boldsymbol a}}
\def\bb{{\boldsymbol b}}
\def\bc{{\boldsymbol c}}
\def\bm{{\boldsymbol m}}
\def\bk{{\boldsymbol k}}
\def\bi{{\boldsymbol i}}
\def\bl{{\boldsymbol l}}
\def\bq{{\boldsymbol q}}
\def\bu{{\boldsymbol u}}
\def\bt{{\boldsymbol t}}
\def\bs{{\boldsymbol s}}
\def\bv{{\boldsymbol v}}
\def\bw{{\boldsymbol w}}
\def\bx{{\boldsymbol x}}
\def\bX{{\boldsymbol X}}
\def\bz{{\boldsymbol z}}
\def\bwy{{\boldsymbol y}}
\def\bY{{\boldsymbol Y}}
\def\bL{{\boldsymbol L}}
\def\baa{{\boldsymbol\alpha}}
\def\bbb{{\boldsymbol\beta}}
\def\bet{{\boldsymbol\eta}}
\def\bxi{{\boldsymbol\xi}}
\def\bo{{\boldsymbol 0}}
\def\bol{{\boldkey 1}_L}
\def\ep{\varepsilon}
\def\p{\boldsymbol\varphi}
\def\q{\boldsymbol\psi}
\def\rank{\operatorname{rank}}
\def\aut{\operatorname{Aut}}
\def\lcm{\operatorname{lcm}}
\def\sgn{\operatorname{sgn}}
\def\spn{\operatorname{span}}
\def\md{\operatorname{mod}}
\def\Norm{\operatorname{Norm}}
\def\dim{\operatorname{dim}}
\def\det{\operatorname{det}}
\def\Vol{\operatorname{Vol}}
\def\rk{\operatorname{rk}}
\def\Gal{\operatorname{Gal}}
\def\WR{\operatorname{WR}}
\def\WO{\operatorname{WO}}
\def\GL{\operatorname{GL}}
\def\pr{\operatorname{pr}}
\def\Tr{\operatorname{Tr}}


\section{Introduction and statement of results}
\label{intro}

Let $K$ be a number field of degree $d = [K:\que] \geq 1$ and $\O_K$ its ring of integers. An element $\theta \in K$ is called {\it primitive} if $K = \que(\theta)$. This is equivalent to the condition that $\deg_{\que}(\theta) = d$, and hence, there are infinitely many primitive elements in $K$. A conjecture of Ruppert \cite{ruppert} (also see \cite{vaaler_widmer-1} for the convenient formulation we are using) asserts that there exists a primitive element $\theta \in K$ such that
\begin{equation}
\label{rupp}
h(\theta) \leq c(d) |\Delta_K|^{\frac{1}{2d}},
\end{equation}
where $h$ is the absolute Weil height, $\Delta_K$ is the discriminant of the number field $K$, and $c(d)$ is a constant depending only on the degree $d$; we review all the necessary notation in Section~\ref{tools}. Ruppert himself proved this conjecture for quadratic number fields and for totally real fields of prime degree. There has been quite a bit of later work on this conjecture; for instance, Vaaler and Widmer \cite{vaaler_widmer-1} proved the conjecture for number fields with at least one real embedding (further results in the case of totally complex fields were just recently obtained in~\cite{avw}). More generally, a slightly weaker bound is obtained by Pazuki and Widmer in \cite[Lemma~7.1]{widmer}:
\begin{equation}
\label{wd0}
h(\theta) \leq |\Delta_K|^{\frac{1}{d}},
\end{equation}
where $\theta$ can be taken in $\O_K$. If $\theta$ is a primitive element, then $1,\theta,\dots,\theta^{d-1}$ is a basis for $K$ as a $\que$-vector space. Hence, Ruppert's conjecture implies the existence of such a basis with 
\begin{equation}
\label{ruppert}
\max_{0 \leq k \leq d-1} h(\theta^k) \leq c(d)^{d-1} |\Delta_K|^{\frac{d-1}{2d}}.
\end{equation}

Consider the situation when $K/\que$ is a Galois extension. In that case, there exists a {\it normal basis} for $K$ over $\que$, i.e. a basis consisting of algebraic conjugates $\beta_1,\dots,\beta_d$; in fact, there are infinitely many such bases (this fact is known as the {\it normal basis theorem}, usually attributed to the works of Noether and Deuring, 1932). On the other hand, if $\beta_1$ is just an arbitrary primitive element for $K$, its algebraic conjugates may be linearly dependent, and so not every primitive element gives rise to a normal basis. Indeed, it may happen for instance that the degree-$d$ minimal polynomial of $\beta_1$ has zero coefficient in front of $x^{d-1}$, which implies that
$$\beta_1 + \dots + \beta_d = 0.$$
It is then natural to ask for a normal basis of bounded height. The first simple observation about quadratic fields follows directly from Ruppert's bound.

\begin{prop} \label{qnb} For all but at most finitely many quadratic extensions $K/\que$, there exists a normal basis $\beta_1,\beta_2 \in K$ with
$$h(\beta_i) \leq c(2) |\Delta_K|^{\frac{1}{4}},$$
for $i=1,2$, where $c(2)$ is as in~\eqref{rupp}.
\end{prop}

\noindent
We give a quick proof of this proposition in Section~\ref{tools}, showing that there can be at most finitely many exceptional quadratic fields $K/\que$ with small discriminant for which this result does not hold. Our first main result produces a general bound for any Galois extension $K/\que$.

\begin{thm} \label{main1} Let $K/\que$ be a Galois extension of degree $d \geq 2$. There exists a normal basis $\beta_1,\dots,\beta_d$ for $K$ over $\que$ so that
$$h(\beta_i) \leq  \frac{d^{4d} (d^2-d+2)^{4d-3}}{2^{4d-3}} \binom{d-1}{[(d-1)/2]}^2 |\Delta_K|^{(d-1)(4d-3)},$$
for all $1 \leq i \leq d$.
\end{thm}

\noindent
Our argument for the proof of this theorem follows the standard proof of the {\it normal basis theorem}. We are using the Pazuki and Widmer bound~\eqref{wd0}, a polynomial non-vanishing principle (Lemma~\ref{nonvanish}) along with some standard inequalities on height and Mahler measure to make this argument effective. We present our proof in Section~\ref{pf_main1}. In the case when $d$ is an odd prime, we can obtain a better bound using a completely different approach.

\begin{thm} \label{main2} Let $K/\que$ be a Galois extension of prime degree $d \geq 3$. There exists a normal basis $\beta_1,\dots,\beta_d$ for $K$ over $\que$ consisting of algebraic integers so that
$$h(\beta_i) \leq |\Delta_K|^{1/2},$$
for all $1 \leq i \leq d$.
\end{thm}

\noindent
We prove Theorem~\ref{main2} in Section~\ref{pf_main2}. Our argument is based on a result of Dubickas~\cite{dub} on linear independence of algebraic conjugates of prime degree and Minkowski's {\it successive minima theorem}. 

It is worth mentioning that our approach in the proofs of Theorems~\ref{main1} and~\ref{main2} utilizes the Diophantine avoidance method, i.e., obtaining effective height-bounds for points {\it not} satisfying some algebraic conditions. Such avoidance ideas are naturally embedded in the proofs of some classical theorems in algebraic number theory, such as the {\it primitive element theorem} and the {\it normal basis theorem}. While most points do not satisfy any given polynomial equation, explicitly identifying such a point may require some work; it is the ``searching for hay in a haystack" problem. We are ready to proceed.
\bigskip

\section{Notation and heights}
\label{tools}

Throughout this paper, we work over a Galois number field $K$ of degree $d \geq 2$ over $\que$ and write $\O_K$ for its ring of integers. Let $\sigma_1,\dots,\sigma_d : K \hookrightarrow \cee$ be the embeddings of $K$. Since $K/\que$ is Galois, $K$ is either totally real or totally imaginary, meaning that either all of the embeddings are real or all of them are complex coming in conjugate pairs.

We normalize absolute values and introduce the standard height function. Let us write $M(K)$ for the set of places of $K$. For each $v \in M(K)$ let $d_v = [K_v : \que_v]$ be the local degree, then for each $u \in M(\que)$, $\sum_{v \mid u} d_v = d$. We select the absolute values so that $|\ |_v$ extends the usual archimedean absolute value on $\que$ when $v \mid \infty$, or the usual $p$-adic absolute value on $\que$ when $v \nmid \infty$. Then archimedean places are in bijective correspondence with the embeddings so that for each $v \mid \infty$ there exists an index $1 \leq j \leq d$ with
$$|\alpha|_v = |\sigma_j(\alpha)|,$$
for each $\alpha \in K$, where $|\ |$ is the usual absolute value on $\real$ or $\cee$ (notice that each conjugate pair of complex embeddings induces the same place). With this normalization choice, the product formula reads
$$\prod_{v \in M(K)} |\alpha|_v^{d_v} = 1,$$
for each nonzero $\alpha \in K$. We define the multiplicative Weil height on algebraic numbers $\alpha \in K$ as
$$h(\alpha) = \prod_{v \in M(K)} \max \{1,|\alpha|_v \}^{d_v/d},$$
and notice that $h(\alpha) = \prod_{v \mid \infty} \max \{1,|\alpha|_v \}^{d_v/d}$ if $\alpha \in \O_K$. This height is absolute, meaning that it is the same when computed over any number field $K$ containing~$\alpha$: this is due to the normalizing exponent~$1/d$ in the definition. Hence, we can compute height for elements of $\qbar$.
\smallskip

We review some useful well-known properties of heights. The first can be found, for instance, as Lemma~2.1 of~\cite{lf3}.

\begin{lem} \label{ht_sum} Let $\xi_1,\dots,\xi_m \in \zed$ and $\alpha_1,\dots,\alpha_m \in \qbar$ for $m \geq 1$. Then
$$h \left( \sum_{j=1}^m \xi_j \alpha_j \right) \leq m |\bxi| \prod_{j=1}^m h(\alpha_j),$$
where $\bxi = (\xi_1,\dots,\xi_m)$ and $|\bxi| := \max \{ |\xi_i| : 1 \leq i \leq m \}$.
\end{lem}

Define {\it Mahler measure} of a polynomial $f(x) = \sum_{k=0}^d a_k x^k \in \cee[x]$ of degree $d$ with roots $\alpha_1,\dots,\alpha_d \in \cee$ to be
$$\mu(f) = |a_d| \prod_{k=1}^d \max \{1,|\alpha_k|\}.$$
The next lemma is Proposition~1.6.6 of~\cite{bg}.

\begin{lem} \label{mahler1} Let $\alpha \in \qbar$ have degree $d$ and let $f(x) \in \zed[x]$ be its minimal polynomial. Then
$$\mu(f) = h(\alpha)^d.$$
\end{lem}

\noindent
Further, write $|f| = \max_{0 \leq n \leq d} |a_d|$, then Lemma~1.6.7 of~\cite{bg} provides the following bound.

\begin{lem} \label{mahler2} Write $[\ ]$ for the integer part function, then
$$|f| \leq \binom{d}{[d/2]} \mu(f),$$
for any $f(x) \in \cee[x]$.
\end{lem}

The next lemma quantifies the basic principle that a polynomial which is not identically zero cannot vanish ``too much". Somewhat different formulations of this principle can be found in~\cite{cassels} (Lemma~1 on p.~261) as well as in the context of N. Alon's celebrated Combinatorial Nullstellensatz~\cite{alon}. The following formulation, which is most convenient for our purposes follows easily from Lemma~2.2 of~\cite{lf1}.

\begin{lem} \label{nonvanish} Let $K$ be a number field as above and $P(\bx) \in K[x_1,\dots,x_n]$ be a polynomial of degree $m$ in $n$ variables which is not identically $0$. There exists a point~$\bxi \in \zed^n$ such that $P(\bxi) \neq 0$ and 
$$|\bxi| \leq \frac{m+2}{2}.$$
\end{lem}

In case $K$ is totally real, we define the Minkowski embedding $\Sigma_K = (\sigma_1,\dots,\sigma_d) : K \hookrightarrow \real^d$, then for any ideal $I \subseteq \O_K$ the image $\Sigma_K(I)$ is a lattice of full rank in $\real^d$. We define the determinant of a full-rank lattice to be the absolute value of the determinant of any basis matrix for the lattice, then
\begin{equation}
\label{ideal_det}
\det(\Sigma_K(I)) = \Nn_K(I) |\Delta_K|^{1/2},
\end{equation}
where $\Nn_K(I) := |\O_K/I|$ is the norm of $I$, as follows, for instance, from Corollary 2.4 of \cite{bayer}. The following property is Lemma 4.1 of~\cite{lf_siki}.

\begin{lem} \label{ht_sigma} Let $K$ be a totally real number field. For any nonzero $\alpha \in \O_K$,
$$1 \leq h(\alpha) \leq |\Sigma_K(\alpha)|,$$
where $|\ |$ stands for the sup-norm on $\real^d$, as above.
\end{lem}

\noindent
We are only using the Minkowski embedding, ideal lattice construction and Lemma \ref{ht_sigma} in Section~\ref{pf_main2} where the $[K:\que]$ is an odd prime, which implies that $K$ is a totally real cyclic extension of $\que$. This is the reason why we are only introducing this notation in the totally real case, where it is simpler than in general.
\smallskip

We finish this section with a proof of Ruppert's bound on the height of a normal basis in the quadratic case.

\proof[Proof of Proposition~\ref{qnb}] Let $K = \que(\sqrt{D})$ for a nonzero squarefree integer $D \neq 1$. A result of Ruppert~\cite{ruppert} guarantees the existence of a primitive element $\theta \in K$ satisfying~\eqref{rupp}, i.e.
$$h(\theta) \leq c(2) |\Delta_K|^{\frac{1}{4}},$$
for an absolute constant $c(2)$. Then $\theta = a+b\sqrt{D}$ for some $a,b \in \que$. Suppose that $a=0$ and $b=m/n$ for some relatively prime integers $m,n$, then the minimal polynomial of $\theta$ is $f(x) = n^2 x^2 - Dm^2$ and, by Lemma~\ref{mahler1}, we have
$$h(\theta) = \mu(f)^{1/2} = |n| \max \{1, |m| |\sqrt{D}| / |n| \} \geq |\sqrt{D}| = c |\Delta_K|^{1/2},$$
for an absolute constant $c$. This implies that we must have $a \neq 0$, in which case $a \pm b\sqrt{D}$ is the desired normal basis, unless $|\Delta_K|^{1/4} \leq c(2)/c$. Hence, there can be at most finitely many exceptions.
\endproof

\bigskip

\section{Proof of Theorem~\ref{main1}}
\label{pf_main1}

We follow the standard proof of the {\it normal basis theorem}, e.g.~\cite[Theorem~28]{artin}, making it effective. Let $G = \{ \sigma_1,\dots,\sigma_d \}$ be the Galois group of $K/\que$ with $\sigma_1$ being the identity, where we are identifying elements of $G$ with the embeddings of $K$ into $\cee$. Let $\alpha \in \O_K$ be a primitive element, $f(x) \in \zed[x]$ be the minimal polynomial of $\alpha$ and define
$$g(x) = \frac{f(x)}{(x-\alpha) f'(\alpha)}.$$
Let $D(x) = \det \left( \sigma_i \sigma_j (g(x)) \right)$. This is a nonzero polynomial with integer coefficients and
$$\deg (D(x)) = d(d-1).$$
We want to choose $\alpha \in \O_K$ so that $D(\alpha) \neq 0$: if this is the case, then the conjugates of $g(\alpha)$ are $\que$-linearly independent, hence, form a normal basis for $K$ over $\que$. We proceed as follows. Let $\theta \in \O_K$ be a primitive element satisfying~\eqref{wd0}, then $1,\theta,\dots,\theta^{d-1} \in \O_K$ is a basis for $K$ over $\que$ with
\begin{equation}
\label{pfm1}
\max_{0 \leq k \leq d-1} h(\theta^k) \leq |\Delta_K|^{\frac{d-1}{d}}.
\end{equation}
For a given vector $\bxi = (\xi_0,\dots,\xi_{d-1}) \in \zed^d$, define
\begin{equation}
\label{pfm2}
\alpha_{\bxi} = \sum_{k=0}^{d-1} \xi_k \theta^k \in \O_K.
\end{equation}
Then $D(\alpha_{\bxi})$ is a polynomial in $d$ variables $\xi_0,\dots,\xi_{d-1}$ of degree $d(d-1)$, which is not identically zero. Then Lemma~\ref{nonvanish} guarantees the existence of an integer vector $\bxi \in \zed^d$ so that $D(\alpha_{\bxi}) \neq 0$ with
$$|\bxi| \leq \frac{d(d-1)+2}{2}.$$
Combining this observation with Lemma~\ref{ht_sum}, \eqref{pfm1} and~\eqref{pfm2}, we obtain an element $\alpha_{\bxi} \in K$ so that $D(\alpha_{\bxi}) \neq 0$ with
\begin{equation}
\label{alpha_bnd}
h(\alpha_{\bxi}) \leq \frac{d(d^2-d+2)}{2} |\Delta_K|^{d-1}.
\end{equation}
For this choice of $\alpha_{\bxi}$, define $\beta = g(\alpha_{\bxi})$. Then, by the argument in the proof of~\cite[Theorem~28]{artin},
$$\beta_j = \sigma_j(\beta),\ 1 \leq j \leq d$$
is a normal basis for $K$. 

We now want to estimate the height of $\beta$. Let us write $\alpha^j_{\bxi} = \sigma_j(\alpha_{\bxi})$ and notice that for each $1 \leq j \leq d$,
$$g_j(x) = \sigma_j(g(x)) = \frac{f(x)}{(x-\alpha^j_{\bxi}) f'(\alpha^j_{\bxi})}$$
has degree $d-1$, roots $\alpha^i_{\bxi}$ for all $i \neq j$, and leading coefficient $1/f'(\alpha^j_{\bxi})$. Then, by Lemma~\ref{mahler2},
$$|g_j| \leq \binom{d-1}{[(d-1)/2]} \mu(g_j) = \binom{d-1}{[(d-1)/2]} \frac{h(\alpha_{\bxi})^d}{|f'(\alpha^j_{\bxi})| \max \{1,|\alpha^j_{\bxi}|\}}.$$
Then, for each $1 \leq j \leq d$,
$$|\beta_j| \leq d |g_j| \max \{ 1, |\alpha^j_{\bxi}| \}^{d-1} \leq d \binom{d-1}{[(d-1)/2]} \frac{h(\alpha_{\bxi})^d \max \{ 1, |\alpha^j_{\bxi}| \}^{d-2}}{|f'(\alpha^j_{\bxi})|}.$$
This implies that
\begin{eqnarray}
\label{beta_1}
\max \{1, |\beta_j| \} & \leq & \frac{1}{|f'(\alpha^j_{\bxi})|} \max \left\{ |f'(\alpha^j_{\bxi})|, d \binom{d-1}{[(d-1)/2]} h(\alpha_{\bxi})^d \max \{ 1, |\alpha^j_{\bxi}| \}^{d-2} \right\} \nonumber \\
& \leq & d \binom{d-1}{[(d-1)/2]} h(\alpha_{\bxi})^d \frac{1}{|f'(\alpha^j_{\bxi})|} \max \{ 1, |f'(\alpha^j_{\bxi})| \} \max \{ 1, |\alpha^j_{\bxi}| \}^{d-2}.
\end{eqnarray}
Notice that the coefficients of $f'(\alpha_{\bxi}) g(x)$, while not necessarily rational integers, are algebraic integers, as is $\alpha_{\bxi}$. Therefore, for any $v \nmid \infty$, 
$$|f'(\alpha_{\bxi})|_v |g|_v \leq 1,$$
where we write $|g|_v$ for the maximum of absolute values $|\ |_v$ of the coefficients of~$g(x)$. This implies that for every $v \nmid \infty$,
\begin{eqnarray}
\label{beta_2}
\max \{ 1, |\beta|_v \} & \leq & \max \{ 1, |g|_v \} \max \{ 1, |\alpha_{\bxi}|_v \}^{d-1} \nonumber \\ 
& \leq & \frac{1}{|f'(\alpha_{\bxi})|_v} \max \{ 1, |f'(\alpha_{\bxi})|_v \} \max \{ 1, |\alpha_{\bxi}|_v \}^{d-1}.
\end{eqnarray}
Observe that $\max \{ 1, |\alpha_{\bxi}|_v \} = 1$ for $v \nmid \infty$, since $\alpha_{\bxi} \in \O_K$. Now, we combine~\eqref{beta_1} and~\eqref{beta_2} and take a product. Using the product formula, we obtain a bound
\begin{equation}
\label{beta3}
h(\beta) = \left\{ \prod_{v \in M(K)} \max \{ 1, |\beta|_v \}^{d_v} \right\}^{\frac{1}{d}} \leq d \binom{d-1}{[(d-1)/2]} h(\alpha_{\bxi})^d h(\alpha_{\bxi})^{d-2} h(f'(\alpha_{\bxi})).
\end{equation}
Notice that $f'(x)$ is a polynomial of degree $d-1$ with integer coefficients and
$$|f'| \leq d |f| \leq d \binom{d-1}{[(d-1)/2]} \mu(f) = d \binom{d-1}{[(d-1)/2]} h(\alpha_{\bxi})^d,$$
by Lemmas~\ref{mahler1} and~\ref{mahler2}, since $f(x)$ is the minimal polynomial of $\alpha_{\bxi}$. This inequality implies that
\begin{eqnarray}
\label{f-prime}
h(f'(\alpha_{\bxi})) & = & \left\{ \prod_{v \mid \infty} \max \{ 1,|f'(\alpha_{\bxi})|_v \}^{d_v} \times \prod_{v \nmid \infty} \max \{ 1,|f'(\alpha_{\bxi})|_v \}^{d_v} \right\}^{\frac{1}{d}} \nonumber \\
& \leq & d |f'| h(\alpha_{\bxi})^{d-1} \leq d^2 \binom{d-1}{[(d-1)/2]} h(\alpha_{\bxi})^{2d-1}.
\end{eqnarray}
Combining~\eqref{beta3} with~\eqref{f-prime} and~\eqref{alpha_bnd}, we obtain
\begin{eqnarray}
\label{fin_bnd}
h(\beta) & \leq & d^3 \binom{d-1}{[(d-1)/2]}^2 h(\alpha_{\bxi})^{4d-3} \nonumber \\
& \leq & \frac{d^{4d} (d^2-d+2)^{4d-3}}{2^{4d-3}} \binom{d-1}{[(d-1)/2]}^2 |\Delta_K|^{(d-1)(4d-3)}.
\end{eqnarray}
Since $h(\beta_j) = h(\beta)$ for every $1 \leq j \leq d$, this completes the proof.
\medskip

\noindent
\begin{rem} \label{better_bnd} Notice that in our proof of Theorem~\ref{main1}, we used the weaker bound~\eqref{wd0} instead of the conjectured bound~\eqref{rupp}, which has been established in the case of real Galois fields. Using that stronger bound would slightly improve the constant depending on $d$ and divide the exponent on $|\Delta_K|$ by~$2$ in the inequality~\eqref{fin_bnd}.
\end{rem}

\bigskip

\section{Proof of Theorem~\ref{main2}}
\label{pf_main2}

Since we are assuming that $K/\que$ is a Galois extension of prime degree $d \geq 3$, then the order of the Galois group $G$ is an odd prime, so $G$ is cyclic; hence, $K/\que$ is a cyclic extension. Further, $K$ must be a totally real number field and so $d_v = 1$ for every archimedean place~$v$. Therefore, the Minkowski embedding notation we introduced in Section~\ref{tools} applies here. Let $L_K = \Sigma_K(\O_K)$ be the full-rank lattice in $\real^d$, then each element $\baa = (\alpha_1,\dots,\alpha_d) \in L_K$ where the coordinates $\alpha_1,\dots,\alpha_d$ are algebraic conjugates. Since $\alpha_i$'s are algebraic integers, we have $|\alpha_i|_v \leq 1$ for every $1 \leq i \leq d$ and $v \nmid \infty$, hence, by product formula, for each $\baa \in L_K$,
\begin{eqnarray}
\label{max0}
|\baa| & = & \max \{|\alpha_1|,\dots,|\alpha_d|\} \geq \left( \prod_{i=1}^d |\alpha_i| \right)^{\frac{1}{d}} = \left( \prod_{v \mid \infty} |\alpha_1|_v \right)^{\frac{1}{d}} \nonumber \\
& \geq & \left( \prod_{v \in M(K)} |\alpha_1|_v^{d_v} \right)^{\frac{1}{d}} = 1.
\end{eqnarray}
We want to find an element $\baa \in L_K$ whose coordinates are linearly independent over~$\que$. We will use the following special case of a result of A. Dubickas, which we state specifically over~$\que$.

\begin{thm} [\cite{dub}, Theorem~1] \label{a_dub} Let $d$ be prime and $c_1,\dots,c_d \in \que$. Then
$$c_1 \alpha_1 + \dots + c_d \alpha_d \in \que$$
if and only if $c_1 = \dots = c_d$.
\end{thm}

\noindent
Hence, this theorem implies that if $\baa = (\alpha_1,\dots,\alpha_d) \in L_K$ is such that
\begin{equation}
\label{lin_f}
\alpha_1 + \dots + \alpha_d \neq 0,
\end{equation}
then $\alpha_1,\dots,\alpha_d \in \O_K$ are linearly independent over~$\que$, thus form a normal basis for $K/\que$. Let
$$\C = \left\{ \bx \in \real^d : |\bx| \leq 1 \right\}$$
be the cube of side-length $2$ centered at the origin in $\real^d$, then $\Vol_d(\C) = 2^d$. Write $\lambda_1,\dots,\lambda_d$ for the successive minima of $\C$ with respect to our lattice $L_K$, then Minkowski's {\it successive minima theorem} provides the bound
\begin{equation}
\label{mink0}
\prod_{i=1}^d \lambda_i \leq \frac{2^d \det(L_K)}{\Vol_d(\C)} = \sqrt{|\Delta_K|},
\end{equation}
by~\eqref{ideal_det}, and
\begin{equation}
\label{max1}
1 \leq \lambda_1 \leq \dots \leq \lambda_d,
\end{equation}
by~\eqref{max0}. Let $\baa_1,\dots,\baa_d \in L_K$ be the linearly independent points corresponding to $\lambda_1,\dots,\lambda_d$, then at least one of these vectors satisfies condition~\eqref{lin_f}, since
$$\left\{ \baa = (\alpha_1,\dots,\alpha_d) \in L_K : \alpha_1 + \dots + \alpha_d = 0 \right\}$$
is a sublattice of rank $d-1$. Let us write $\bbb = (\beta_1,\dots,\beta_d)$ for the $\baa_i$ satisfying~\eqref{lin_f}, then $\beta_1,\dots,\beta_d \in \O_K$ is a normal basis for $K/\que$ and~\eqref{mink0} combined with~\eqref{max1} imply that
\begin{equation}
\label{mink1}
\max \{ |\beta_1|,\dots,|\beta_d| \} \leq |\Delta_K|^{1/2}.
\end{equation}
Combining~\eqref{mink1} with Lemma~\ref{ht_sigma}, we obtain
$$h(\beta_i) \leq |\Delta_K|^{1/2},$$
for each $1 \leq i \leq d$. This completes the proof.
\bigskip

\noindent
{\bf Acknowledgements:} We wish to thank Professor Adebisi Agboola for a valuable suggestion that led us to consider an effective version of the {\it normal basis theorem}. We also thank the referee for a thorough review of our paper and helpful suggestions.

\bigskip

\noindent
{\bf Data availability statement:} Data sharing not applicable to this article as no datasets were generated or analyzed during the current study.
\bigskip

\noindent
{\bf Conflict of interest statement:} The authors declare no conflict of interest.
\bigskip

\bibliographystyle{plain}  

\end{document}